\input amstex
\input amsppt.sty
\magnification=\magstep1
\hsize=32truecc
\vsize=22.2truecm
\baselineskip=16truept
\NoBlackBoxes
\TagsOnRight \pageno=1 \nologo
\def\Z{\Bbb Z}
\def\N{\Bbb N}

\def\l{\left}
\def\r{\right}
\def\bg{\bigg}
\def\({\bg(}
\def\[{\bg\lfloor}
\def\){\bg)}
\def\]{\bg\rfloor}
\def\t{\text}
\def\f{\frac}

\def\bi{\binom}
\def\eq{\equiv}

\def\ls{\leqslant}
\def\gs{\geqslant}
\def\mo{\roman{mod}}

\def\M#1#2{\thickfracwithdelims[]\thickness0{#1}{#2}_q}

\def\Proof{\noindent{\it Proof}}

\def\Remark{\medskip\noindent{\it  Remark}}

\def\Ack{\medskip\noindent {\bf Acknowledgment}}
\hbox {Ramanujan J. 40(2016), no.\,3, 511-533.}
\bigskip
\topmatter
\title Congruences involving $g_n(x)=\sum_{k=0}^n\bi nk^2\bi{2k}kx^k$\endtitle
\author Zhi-Wei Sun\endauthor
\leftheadtext{Zhi-Wei Sun}
\affil Department of Mathematics, Nanjing University\\
 Nanjing 210093, People's Republic of China
  \\  zwsun\@nju.edu.cn
  \\ {\tt http://math.nju.edu.cn/$\sim$zwsun}
\endaffil
\abstract Define $g_n(x)=\sum_{k=0}^n\binom nk^2\binom{2k}kx^k$
 for $n=0,1,2,\ldots$.  Those numbers $g_n=g_n(1)$ are closely related to Ap\'ery numbers and Franel numbers.
 In this paper we establish some fundamental congruences
 involving $g_n(x)$. For example, for any prime $p>5$ we have
 $$\sum_{k=1}^{p-1}\f{g_k(-1)}k\eq0\pmod{p^2}\quad\text{and}\quad\sum_{k=1}^{p-1}\f{g_k(-1)}{k^2}\eq0\pmod p.$$
 This is similar to Wolstenholme's classical congruences
 $$\sum_{k=1}^{p-1}\f1k\eq0\pmod{p^2}\quad\text{and}\quad\sum_{k=1}^{p-1}\f{1}{k^2}\eq0\pmod p$$
 for any prime $p>3$.
\endabstract
\thanks 2010 {\it Mathematics Subject Classification}. \,Primary 11A07, 11B65;
Secondary  05A10, 05A30, 11B75.
\newline\indent {\it Keywords}. Franel numbers, Ap\'ery numbers,
binomial coefficients, congruences.
\newline\indent Supported by the National Natural Science
Foundation (Grant No. 11571162) of China.
\endthanks

\endtopmatter
\document

\heading{1. Introduction}\endheading

It is well known that
$$\sum_{k=0}^n\bi nk^2=\bi{2n}n\ \ (n=0,1,2,\ldots)$$
and central binomial coefficients play important roles in
mathematics. A famous theorem of J. Wolstenholme [W] asserts that for any prime $p>3$ we have
$$\f12\bi{2p}p=\bi{2p-1}{p-1}\eq1\pmod{p^3},$$
$$ H_{p-1}\eq0\pmod{p^2}\quad\t{and}\quad H_{p-1}^{(2)}\eq0\pmod p,$$
where
$$H_n:=\sum_{0<k\ls n}\f1k\quad\t{and}\quad H_n^{(2)}:=\sum_{0<k\ls n}\f1{k^2}\quad\t{for}\ n\in\N=\{0,1,2,\ldots\};$$
see also [Zh] for some extensions.
The reader may consult [S11a], [S11b], [ST1] and [ST2] for recent work on congruences involving central binomial coefficients.

The Franel numbers given by
$$f_n=\sum_{k=0}^n\bi nk^3\ \ \ (n=0,1,2,\ldots)$$
(cf. [Sl, A000172]) were first introduced by J. Franel in 1895 who noted the recurrence relation:
$$(n+1)^2f_{n+1}=(7n(n+1)+2)f_n+8n^2f_{n-1}\ (n=1,2,3,\ldots).$$
In 1992 C. Strehl [St92] showed that the Ap\'ery numbers given by
$$A_n=\sum_{k=0}^n\bi nk^2\bi{n+k}k^2=\sum_{k=0}^n\bi{n+k}{2k}^2\bi{2k}k^2\
(n=0,1,2,\ldots)$$ (arising from Ap\'ery's proof of the irrationality of $\zeta(3)=\sum_{n=1}^\infty1/n^3$ (cf. [vP])) can be expressed in terms of Franel numbers, namely,
$$A_n=\sum_{k=0}^n\bi nk\bi{n+k}kf_k.\tag1.1$$

Define
$$g_n=\sum_{k=0}^n\bi nk^2\bi{2k}k\quad\t{for}\ n\in\N.\tag1.2$$
Such numbers are interesting due to Barrucand's identity ([B])
$$\sum_{k=0}^n\bi nkf_k=g_n\quad (n=0,1,2,\ldots).\tag1.3$$
For a combinatorial interpretation of such numbers, see D. Callan [C].
The sequences $(f_n)_{n\gs0}$ and $(g_n)_{n\gs0}$ are two of the five sporadic sequences (cf. D. Zagier [Z, Section 4])
which are integral solutions of certain Ap\'ery-like recurrence equations and closely related to the theory of modular forms.

In [S12] and [S13b] the author introduced the Ap\'ery polynomials
$$A_n(x):=\sum_{k=0}^n\bi nk^2\bi{n+k}k^2x^k\ \ (n=0,1,2,\ldots)$$
and the Franel polynomials
$$f_n(x):=\sum_{k=0}^n\bi nk^2\bi{2k}nx^k=\sum_{k=0}^n\bi nk\bi k{n-k}\bi{2k}kx^k\ (n=0,1,2,\ldots),$$
and deduced various congruences involving such polynomials.
(Note that $A_n(1)=A_n$, and $f_n(1)=f_n$ by [St94].) See also [S13a] for connections between primes $p=x^2+3y^2$ and
the Franel numbers. Here we introduce the polynomials
$$g_n(x):=\sum_{k=0}^n\bi nk^2\bi{2k}kx^k\ \ (n=0,1,2,\ldots).$$
Both $f_n(x)$ and $g_n(x)$ play important roles in some kinds of series for $1/\pi$ (cf. Conjecture 3 and the subsequent remark in [S11]).

In this paper we study various congruences involving $g_n(x)$. As usual, for an odd prime $p$ and an integer $a$,
$(\f ap)$ denotes the Legendre symbol, and $q_p(a)$ stands for the Fermat quotient $(a^{p-1}-1)/p$
if $p\nmid a$. Also, $B_0,B_1,B_2,\ldots$ are the well-known Bernoulli numbers and $E_0,E_1,E_2,\ldots$ are the Euler numbers.

Now we state our main results.

\proclaim{Theorem 1.1} Let $p>3$ be a prime.

{\rm (i)} We have
$$\sum_{k=0}^{p-1}g_k(x)(1-p^2H_k^{(2)})\eq \sum_{k=0}^{p-1}\f{p}{2k+1}\l(1-2p^2H_k^{(2)}\r)x^k\pmod{p^4}.\tag1.4$$
Consequently,
$$\align\sum_{k=1}^{p-1}g_k\eq&p^2\sum_{k=1}^{p-1}g_kH_k^{(2)}+\f 76p^3B_{p-3}\pmod{p^4},\tag1.5
\\\sum_{k=0}^{p-1}g_k(-1)\eq&\l(\f{-1}p\r)+p^2\(\sum_{k=0}^{p-1}g_k(-1)H_k^{(2)}-E_{p-3}\)\ (\mo\ p^3),\tag1.6
\\\sum_{k=0}^{p-1}g_k(-3)\eq&\l(\f p3\r)\pmod{p^2}.\tag1.7
\endalign$$

{\rm (ii)} We also have
$$\align
\sum_{k=1}^{p-1}\f{g_k(x)}k\eq&0\pmod{p},\tag1.8
\\\sum_{k=1}^{p-1}\f{g_{k-1}}k\eq&-\l(\f p3\r)2q_p(3)\pmod p,\tag1.9
\\\sum_{k=1}^{p-1}kg_k\eq&-\f 34\pmod{p^2},\tag1.10
\endalign$$
and moreover
$$\f1{3n^2}\sum_{k=0}^{n-1}(4k+3)g_k=\sum_{k=0}^{n-1}\bi{n-1}k^2C_k\tag1.11$$
for all $n\in\Z^+=\{1,2,3,\ldots\}$, where $C_k$ denotes the Catalan number $\bi{2k}k/(k+1)=\bi{2k}k-\bi{2k}{k+1}$.

{\rm (iii)} Provided $p>5$, we have
$$\align\sum_{k=1}^{p-1}\f{g_k(-1)}{k^2}\eq&0\pmod{p},\tag1.12
\\\sum_{k=1}^{p-1}\f{g_k(-1)}{k}\eq&0\pmod {p^2},\tag1.13
\\\sum_{k=1}^{p-1}\f{(-1)^kf_k(-1)}{k}H_k\eq&-2\l(\f{-1}p\r)E_{p-3}\pmod{p}.\tag1.14
\endalign$$
\endproclaim
\Remark\ 1.1. Let $p>3$ be a prime. By [JV, Lemma 2.7], $g_k\eq(\f p3)9^kg_{p-1-k}\pmod p$ for all $k=0,\ldots,p-1$.
So (1.9) implies that
$$\sum_{k=1}^{p-1}\f{g_k}{k9^k}\eq\l(\f p3\r)\sum_{k=1}^{p-1}\f{g_{p-1-k}}k=\l(\f p3\r)\sum_{k=1}^{p-1}\f{g_{k-1}}{p-k}\eq2q_p(3)\pmod p.$$
We conjecture further that
$$\sum_{k=1}^{p-1}\f{g_{k-1}}k\eq-\l(\f p3\r)q_p(9)\pmod{p^2}\ \ \t{and}\ \ \sum_{k=0}^{p-1}\f{g_k}{9^k}\eq\l(\f p3\r)\pmod{p^2}.$$
In [S13b] the author showed the following congruences similar to (1.12) and (1.13):
$$\sum_{k=1}^{p-1}\f{(-1)^kf_k}{k^2}\eq0\pmod p\quad\t{and}\quad\sum_{k=1}^{p-1}\f{(-1)^kf_k}{k}\eq0\pmod{p^2}.$$
Such congruences are interesting in view of Wolstenholme's congruences $H_{p-1}\eq0\pmod{p^2}$
and $H_{p-1}^{(2)}\eq0\pmod p$. Applying the Zeilberger algorithm (cf. [PWZ, pp.\,101-119]) via {\tt Mathematica 9} we find the recurrence for $s_n=g_n(-1)\ (n=0,1,2,\ldots)$:
$$\align &(n+3)^2(4n+5)s_{n+3}+ (20n^3+125n^2+254n+165)s_{n+2}
\\&\quad +(76n^3+399n^2+678n+375)s_{n+1}-25(n+1)^2(4n+9)s_n = 0.
\endalign$$
In contrast with (1.11), we are also able to show the congruence
$$\sum_{k=0}^{p-1}(3k+1)\f{f_k}{8^k}\eq p^2-2p^3q_p(2)+4p^4q_p(2)^2\pmod{p^5}\tag1.15$$
via the combinatorial identity
$$\f1{n^2}\sum_{k=0}^{n-1}(3k+1)f_k8^{n-1-k}=\sum_{k=0}^{n-1}\bi{n-1}k^3\l(1-\f n{k+1}+\f{n^2}{(k+1)^2}\r)\tag1.16$$
which can be shown by the Zeilberger algorithm.
\medskip

We are going to investigate in the next section connections among the polynomials $A_n(x)$, $f_n(x)$ and $g_n(x)$.
Section 3 is devoted to our proof of Theorem 1.1.
In Section 4 we shall propose some conjectures for further research.

\heading{2. Relations among $A_n(x),f_n(x)$ and $g_n(x)$}\endheading

 Obviously,

 $$\f1{n}\sum_{k=0}^{n-1}(2k+1)=n\in\Z\ \ \t{and}\ \ \f1n\sum_{k=0}^{n-1}(2k+1)(-1)^k=(-1)^{n-1}\in\Z$$
 for all $n=1,2,3,\ldots$. This is a special case of our following general result.

\proclaim{Theorem 2.1} Let
$$X_n=\sum_{k=0}^n\bi nk\bi{n+k}kx_k\quad\t{and}\quad y_n=\sum_{k=0}^n\bi nkx_k\quad\t{for all}\ n\in\N.\tag2.1$$
Then
$$X_n=\sum_{k=0}^n\bi nk\bi{n+k}k(-1)^{n-k}y_k\quad\t{for every}\ n\in\N.\tag2.2$$
Also, for any $n\in\Z^+$ we have
$$\f{(-1)^{n-1}}n\sum_{k=0}^{n-1}(2k+1)X_k=\sum_{k=0}^{n-1}\bi{n-1}k\bi{n+k}k(-1)^ky_k \tag2.3$$
and
$$\f{(-1)^{n-1}}n\sum_{k=0}^{n-1}(2k+1)(-1)^kX_k=\sum_{k=0}^{n-1}\bi{n-1}k\bi{n+k}kx_k.\tag2.4$$
\endproclaim
\Proof.  If $n\in\N$, then
$$\align&\sum_{l=0}^n\bi nl\bi{n+l}l(-1)^ly_l
\\=&\sum_{l=0}^n\bi nl\bi{-n-1}l\sum_{k=0}^l\bi lkx_k
\\=&\sum_{k=0}^n\bi nkx_k\sum_{l=k}^n\bi{n-k}{n-l}\bi{-n-1}l
\\=&\sum_{k=0}^n\bi nkx_k\bi{-k-1}n\ \ (\t{by the Chu-Vandermonde identity [G, (2.1)]})
\\=&(-1)^n\sum_{k=0}^n\bi nk\bi{n+k}kx_k
\endalign$$
and hence (2.2) holds.

For any given integer $k\gs0$, by induction on $n$ we have
$$\sum_{l=k}^{n-1}(-1)^l(2l+1)\bi{l+k}{2k}=(-1)^{n-1}(n-k)\bi{n+k}{2k}\tag2.5$$
for all $n=k+1,k+2,\ldots$. Fix a positive integer $n$. In view of (2.2) and (2.5),
$$\align \sum_{l=0}^{n-1}(2l+1)X_l=&\sum_{l=0}^{n-1}(2l+1)\sum_{k=0}^l\bi{l+k}{2k}\bi{2k}k(-1)^{l-k}y_k
\\=&\sum_{k=0}^{n-1}\bi{2k}k(-1)^ky_k\sum_{l=k}^{n-1}(-1)^l(2l+1)\bi{l+k}{2k}
\\=&\sum_{k=0}^{n-1}\bi{2k}k(-1)^ky_k(-1)^{n-1}(n-k)\bi{n+k}{2k}
\\=&(-1)^{n-1}n\sum_{k=0}^{n-1}\bi{n-1}k\bi{n+k}k(-1)^ky_k.
\endalign$$
This proves (2.3). Similarly,
$$\align \sum_{l=0}^{n-1}(2l+1)(-1)^lX_l=&\sum_{l=0}^{n-1}(2l+1)(-1)^l\sum_{k=0}^l\bi{l+k}{2k}\bi{2k}kx_k
\\=&\sum_{k=0}^{n-1}\bi{2k}kx_k\sum_{l=k}^{n-1}(-1)^l(2l+1)\bi{l+k}{2k}
\\=&\sum_{k=0}^{n-1}\bi{2k}kx_k(-1)^{n-1}(n-k)\bi{n+k}{2k}
\\=&(-1)^{n-1}n\sum_{k=0}^{n-1}\bi{n-1}k\bi{n+k}kx_k.
\endalign$$
and hence (2.4) is also valid.

Combining the above, we have completed the proof of Theorem 2.1. \qed

\proclaim{Lemma 2.1} For any nonnegative integers $m$ and $n$ we have the combinatorial identity
$$\sum_{k=0}^n\bi{m-x+y}k\bi{n+x-y}{n-k}\bi{x+k}{m+n}=\bi xm\bi yn.\tag 2.6$$
\endproclaim
\Remark\ 2.1. (2.6) is  due to Nanjundiah, see, e.g., (4.17) of [G, p.\,53].
\medskip

The author [S12] proved that $\f1n\sum_{k=0}^{n-1}(2k+1)A_k(x)\in\Z[x]$ for all $n\in\Z^+$, and conjectured that
$\f1n\sum_{k=0}^{n-1}(2k+1)(-1)^kA_k(x)\in\Z[x]$ for any $n\in\Z^+$, which was confirmed by Guo and Zeng [GZ].

\proclaim{Theorem 2.2} Let $n$ be any nonnegative integer. Then
$$\sum_{k=0}^n \bi nkf_k(x)=g_n(x),\ \ f_n(x)=\sum_{k=0}^n\bi nk(-1)^{n-k}g_k(x),\tag2.7$$
and
$$A_n(x)=\sum_{k=0}^n\bi nk\bi{n+k}kf_k(x)=\sum_{k=0}^n\bi nk\bi{n+k}k(-1)^{n-k}g_k(x).\tag2.8$$
Also, for any $n\in\Z^+$ we have
$$\f{(-1)^{n-1}}n\sum_{k=0}^{n-1}(2k+1)A_k(x)=\sum_{k=0}^{n-1}\bi{n-1}k\bi{n+k}k(-1)^kg_k(x)\tag2.9$$
and
$$\f{(-1)^{n-1}}n\sum_{k=0}^{n-1}(2k+1)(-1)^kA_k(x)=\sum_{k=0}^{n-1}\bi{n-1}k\bi{n+k}kf_k(x).\tag2.10$$
\endproclaim
\Proof. By the binomial inversion formula (cf. (5.48) of [GKP, p.\,192]), the two identities in (2.7) are equivalent.
Observe that
$$\align\sum_{l=0}^n\bi nlf_l(x)=&\sum_{l=0}^n\bi nl\sum_{k=0}^l\bi lk\bi{k}{l-k}\bi{2k}kx^k
\\=&\sum_{k=0}^n\bi nk\bi{2k}kx^k\sum_{l=k}^n\bi{n-k}{n-l}\bi k{l-k}
\\=&\sum_{k=0}^n\bi nk\bi{2k}kx^k\bi n{n-k}=g_n(x)
\endalign$$
with the help of the Chu-Vandermonde identity. Thus (2.7) holds.

 Next we show (2.8). Clearly
$$\align\sum_{l=0}^n\bi nl\bi{n+l}lf_l(x)=&\sum_{l=0}^n\bi nl\bi{n+l}l\sum_{k=0}^l\bi lk\bi k{l-k}\bi{2k}kx^k
\\=&\sum_{k=0}^n\bi nk\bi{2k}kx^k\sum_{l=k}^n\bi{n-k}{l-k}\bi k{l-k}\bi{n+l}n
\\=&\sum_{k=0}^n\bi nk\bi{2k}kx^k\sum_{j=0}^k\bi{n-k}j\bi k{k-j}\bi{n+k+j}{n}
\\=&\sum_{k=0}^n\bi nk\bi{2k}kx^k\bi{n+k}{n-k}\bi{n+k}k\ \ (\t{by Lemma 2.1}).
\endalign$$
This proves the first identity in (2.8). Applying Theorem 2.1 with $x_n=f_n(x)$ and $X_n=A_n(x)$ for $n\in\N$, we get the identity
$$A_n(x)=\sum_{k=0}^n\bi nk\bi{n+k}k(-1)^{n-k}g_k(x)\tag2.11$$
as well as (2.9) and (2.10), with the help of (2.7).

The proof of Theorem 2.2 is now complete. \qed

\Remark\ 2.2. (2.7) and (2.8) in the case $x=1$ are well known.
\medskip

\proclaim{Corollary 2.1} Let $p$ be an odd prime. Then
$$\sum_{k=0}^{p-1}A_k(x)\eq p\sum_{k=0}^{p-1}\f{(-1)^kf_k(x)}{2k+1}\pmod{p^2}\tag2.12$$
and
$$\sum_{k=0}^{p-1}(-1)^kA_k(x)\eq p\sum_{k=0}^{p-1}\f{g_k(x)}{2k+1}\pmod{p^2}.\tag2.13$$
\endproclaim
\Proof. In view of (2.8),
$$\align\sum_{l=0}^{p-1}A_l(x)=&\sum_{l=0}^{p-1}\sum_{k=0}^l\bi{k+l}{2k}\bi{2k}kf_k(x)
=\sum_{k=0}^{p-1}\bi{2k}kf_k(x)\sum_{l=k}^{p-1}\bi{k+l}{2k}
\\=&\sum_{k=0}^{p-1}\bi{2k}kf_k(x)\bi{p+k}{2k+1}=\sum_{k=0}^{p-1}\bi{2k}kf_k(x)\f p{(2k+1)!}\prod_{0<j\ls k}(p^2-j^2)
\\\eq&\sum_{k=0}^{p-1}f_k(x)\f p{2k+1}(-1)^k\pmod{p^2}.
\endalign$$
Similarly,
$$\align\sum_{l=0}^{p-1}(-1)^lA_l(x)=&\sum_{l=0}^{p-1}\sum_{k=0}^l\bi{k+l}{2k}\bi{2k}k(-1)^kg_k(x)
\\=&\sum_{k=0}^{p-1}\bi{2k}k(-1)^kg_k(x)\bi{p+k}{2k+1}
\\\eq&\sum_{k=0}^{p-1}g_k(x)\f p{2k+1}\pmod{p^2}.
\endalign$$
This concludes the proof of Corollary 2.1. \qed

\Remark\ 2.3. In [S12] the author investigated $\sum_{k=0}^{p-1}(\pm1)^kA_k(x)$ mod $p^2$ (where $p$ is an odd prime)
and made some conjectures.
\medskip

For any $n\in\Z$ we set
$$[n]_q=\f{1-q^n}{1-q}=\cases\sum_{0\ls k<n}q^k&\t{if}\ n\gs0,\\-q^n\sum_{0\ls k<-n}q^k&\t{if}\ n<0;\endcases$$
this is the usual $q$-analogue of the integer $n$. Define
$$\M n0=1\ \ \t{and}\ \ \M nk=\prod_{j=1}^k\f{[n-j+1]_q}{[j]_q}\ \ \t{for}\ k\in\Z^+.$$
Obviously, $\lim_{q\to1}\M nk=\bi nk$.

For $n\in\N$ we define
$$A_n(x;q):=\sum_{k=0}^nq^{2n(n-k)}\M nk^2\M{n+k}k^2x^k$$
and
$$g_n(x;q):=\sum_{k=0}^nq^{2n(n-k)}\M nk^2\M{2k}kx^k.$$
Clearly
$$\lim_{q\to1}A_n(x;q)=A_n(x)\ \ \t{and}\ \ \lim_{q\to1}g_n(x;q)=g_n(x).$$
Those identities in Theorem 2.2 have their $q$-analogues. For example, the following theorem gives a $q$-analogue of (2.11).

\proclaim{Theorem 2.3} Let $n\in\N$. Then we have
$$A_n(x;q)=\sum_{k=0}^n(-1)^{n-k}q^{(n-k)(5n+3k+1)/2}\M nk\M{n+k}kg_k(x;q).\tag2.14$$
\endproclaim
\Proof. Let $j\in\{0,\ldots,n\}$. By the $q$-Chu-Vandermonde identity (see, e.g., Ex. 4(b) of [AAR, p.\,542]),
$$\sum_{k=j}^nq^{(k-j)^2}\M{-n-1-j}{k-j}\M{n-j}{n-k}=\M{-2j-1}{n-j}.$$
This, together with
$$\M{-n-1}k\M kj=\M{-n-1}j\M{-n-1-j}{k-j},$$
yields that
$$\sum_{k=j}^nq^{(k-j)^2}\M{-n-1}k\M kj\M{n-j}{k-j}=\M{-n-1}j\M{-2j-1}{n-j}.$$
It is easy to see that
$$\M{-m-1}k=(-1)^kq^{-km-k(k+1)/2}\M{m+k}k.$$
So we are led to the identity
$$\sum_{k=j}^n (-1)^{n-k}q^{\bi{n-k+1}2+2j(n-k)}\M{n+k}k\M kj\M{n-j}{k-j}=\M{n+j}j\M{n+j}{2j}.\tag2.15$$
Since
$$\M nk\M kj=\M nj\M{n-j}{k-j}\ \ \t{and}\ \ \M{n}j\M{n+j}j=\M{n+j}{2j}\M{2j}j,$$
multiplying both sides of (2.15) by $\M nj\M{2j}jx^j$ we get
$$\sum_{k=j}^n(-1)^{n-k}q^{\bi{n-k+1}2+2j(n-k)}\M nk\M{n+k}k\M kj^2\M{2j}jx^j=\M nj^2\M{n+j}j^2x^j.$$

In view of the last identity we can easily deduce the desired (2.14). \qed

By applying Theorem 2.2 we obtain the following new result.

\proclaim{Theorem 2.4} Let $n$ be any positive integer. Then
$$\sum_{k=0}^{n-1}(-1)^{k}(6k^3+9k^2+5k+1)A_k\eq0\pmod{n^3}.\tag2.16$$
\endproclaim
\Proof. By induction on $n$, for each $k=0,\ldots,n-1$ we have
$$\sum_{l=k}^{n-1}(-1)^l(6l^3+9l^2+5l+1)\bi{l+k}{2k}=(-1)^{n-1}(n-k)(3n^2-3k-2)\bi{n+k}{2k}.$$
Thus, in view of (2.8),
$$\align&\f{1}n\sum_{l=0}^{n-1}(-1)^{n-l}(6l^3+9l^2+5l+1)A_l(x)
\\=&\f{(-1)^n}n\sum_{l=0}^{n-1}(-1)^l(6l^3+9l^2+5l+1)\sum_{k=0}^l\bi{l+k}{2k}\bi{2k}kf_k(x)
\\=&\f{(-1)^n}n\sum_{k=0}^{n-1}\bi{2k}kf_k(x)\sum_{l=k}^{n-1}(-1)^l(6l^3+9l^2+5l+1)\bi{l+k}{2l}
\\=&\f{(-1)^n}n\sum_{k=0}^{n-1}\bi{2k}kf_k(x)(-1)^{n-1}(n-k)(3n^2-3k-2)\bi{n+k}{2k}
\\=&\sum_{k=0}^{n-1}\bi{n-1}k\bi{n+k}k(3k+2-3n^2)f_k(x).
\endalign$$
Hence we have reduced (2.16) to the congruence
$$\sum_{k=0}^{n-1}\bi{n-1}k\bi{n+k}k(3k+2)f_k\eq0\pmod{n^2}.\tag2.17$$

The author [S13a, (1.12)] conjectured that
$$a_m:=\f1{m^2}\sum_{k=0}^{m-1}(3k+2)(-1)^kf_k\in\Z\quad\t{for all}\ m=1,2,3,\ldots,$$
and this was confirmed by V.J.W. Guo [Gu]. Set $a_0=0$. Observe that
$$\align&\sum_{k=0}^{n-1}\bi {n-1}k\bi{n+k}k(3k+2)f_k
\\=&\sum_{k=0}^{n-1}\bi{n-1}k\bi{-n-1}k\l((k+1)^2a_{k+1}-k^2a_k\r)
\\=&\sum_{k=1}^n\bi{n-1}{k-1}\bi{-n-1}{k-1}k^2a_k-\sum_{k=0}^{n-1}\bi{n-1}k\bi{-n-1}kk^2a_k
\\=&\bi{-n-1}{n-1}n^2a_n+\sum_{0<k<n}k^2a_k\l(\bi{n-1}{k-1}\bi{-n-1}{k-1}-\bi{n-1}k\bi{-n-1}k\r).
\endalign$$
As
$$\bi{n-1}{k-1}\bi{-n-1}{k-1}-\bi{n-1}k\bi{-n-1}k=\f{n^2}{k^2}\bi{n-1}{k-1}\bi{-n-1}{k-1}$$
for all $k=1,\ldots,n-1,$
we have (2.17) by the above, and hence (2.16) holds. \qed
\medskip

The author [S12] conjectured that for any prime $p>3$ we have
$$\sum_{k=0}^{p-1}(2k+1)(-1)^kA_k\eq p\l(\f p3\r)\pmod {p^3},\tag 2.18$$
and this was confirmed by Guo and Zeng [GZ].

\proclaim{Corollary 2.2} Let $p>3$ be a prime. Then
$$\sum_{k=0}^{p-1}(2k+1)^3(-1)^kA_k\eq-\f p3\l(\f p3\r)\pmod{p^3}.\tag2.19$$
\endproclaim
\Proof. Clearly
$$3(2k+1)^3=4(6k^3+9k^2+5k+1)-(2k+1).$$
Thus (2.19) follows from (2.16) and (2.18). \qed

\Remark\ 2.4. Let $p>3$ be a prime.
We are also able to prove that
$$\sum_{k=0}^{p-1}(2k+1)^5(-1)^kA_k\eq-\f {13}{27}p\l(\f p3\r)\pmod{p^3}\tag2.20$$
and
$$\sum_{k=0}^{p-1}(2k+1)^7(-1)^kA_k\eq\f 59p\l(\f p3\r)\pmod{p^3}.\tag2.21$$
It seems that for each $r=0,1,2,\ldots$ there is a $p$-adic integer $c_r$ only depending on $r$ such that
$$\sum_{k=0}^{p-1}(2k+1)^{2r+1}(-1)^kA_k\eq c_r p\l(\f p3\r)\pmod{p^3}.$$

\heading{3. Proof of Theorem 1.1}\endheading

\proclaim{Lemma 3.1} For any odd prime $p$, we have
$$\f1p\sum_{k=0}^{p-1}(2k+1)A_k(x)\eq\sum_{k=0}^{p-1}g_k(x)-p^2\sum_{k=0}^{p-1}g_k(x)H_k^{(2)}\pmod{p^4}.\tag3.1$$
\endproclaim
\Proof. Obviously,
$$(-1)^k\bi{p-1}k\bi{p+k}k=\prod_{0<j\ls k}\l(1-\f{p^2}{j^2}\r)\eq1-p^2H_k^{(2)}\pmod{p^4}\tag3.2$$
for every $k=0,\ldots,p-1$. Thus (3.1) follows from (2.9) with $n=p$. \qed

\proclaim{Lemma 3.2} Let $p>3$ be a prime. Then
$$g_{p-1}\eq\l(\f p3\r)(1+2p\,q_p(3))\pmod{p^2}.\tag3.3$$
\endproclaim
\Proof. For $k=0,\ldots,p-1$, clearly
$$\bi{p-1}k^2=\prod_{0<j\ls k}\l(1-\f{p}j\r)^2\eq\prod_{0<j\ls k}\l(1-\f{2p}j\r)
=(-1)^k\bi{2p-1}k\pmod{p^2}.$$
Thus, with the help of  [S12b, Corollary 2.2] we obtain
$$g_{p-1}\eq\sum_{k=0}^{p-1}\bi{2p-1}k(-1)^k\bi{2k}k\eq\l(\f p3\r)\l(2\times 3^{p-1}-1\r)\pmod{p^2}.$$
and hence (3.3) holds. \qed

\proclaim{Lemma 3.3} For any odd prime $p$, we have
$$p\sum_{k=0}^{p-1}\f{(-3)^k}{2k+1}\eq\l(\f p3\r)\pmod{p^2}.\tag3.4$$
\endproclaim
\Proof. Clearly (3.4) holds for $p=3$. Below we assume $p>3$.
Observe that
$$\align\sum^{p-1}\Sb k=0\\k\not=(p-1)/2\endSb\f{(-3)^k}{2k+1}
=&\sum_{k=1}^{(p-1)/2}\(\f{(-3)^{(p-1)/2-k}}{2((p-1)/2-k)+1}+\f{(-3)^{(p-1)/2+k}}{2((p-1)/2+k)+1}\)
\\\eq&\l(\f{-3}p\r)\f12\sum_{k=1}^{(p-1)/2}\l(\f{(-3)^k}k-\f13\cdot\f{(-3)^{p-k}}{p-k}\r)
\\=&\f12\l(\f p3\r)\(\f43\sum_{k=1}^{(p-1)/2}\f{(-3)^k}k-\f13\sum_{k=1}^{p-1}\f{(-3)^k}k\)
\\=&-2\l(\f p3\r)\sum_{k=1}^{(p-1)/2}\f{(-3)^{k-1}}k+\f12\l(\f p3\r)\sum_{k=1}^{p-1}\f{(-3)^{k-1}}k
\pmod{p}.\endalign$$
Since
$$\f1p\bi pk=\f1k\bi{p-1}{k-1}\eq\f{(-1)^{k-1}}k\pmod{p}\quad\t{for}\ k=1,\ldots,p-1,$$
we have
$$\align\sum_{k=1}^{p-1}\f{(-3)^{k-1}}k\eq&\f1{3p}\sum_{k=1}^{p-1}\bi pk3^k=\f{4^p-1-3^p}{3p}= 4(2^{p-1}+1)\f{2^{p-1}-1}{3p}-\f{3^{p-1}-1}p
\\\eq&\f 83q_p(2)-q_p(3)\pmod p.
\endalign$$
Note also that
$$\align\sum_{k=1}^{(p-1)/2}\f{(-3)^{k-1}}k=&\sum_{k=1}^{(p-1)/2}\int_0^1(-3x)^{k-1}dx=\int_0^1\f{1-(-3x)^{(p-1)/2}}{1+3x}dx
\\=&\int_0^1\sum_{k=1}^{(p-1)/2}\bi{(p-1)/2}k(-1-3x)^{k-1}dx
\\=&\sum_{k=1}^{p-1}\bi{(p-1)/2}k\f{(-1-3x)^k}{-3k}\bigg|_{x=0}^1
\\\eq&\sum_{k=1}^{p-1}\bi{-1/2}k\f{(-1)^k-(-4)^k}{3k}
=\f13\sum_{k=1}^{p-1}\f{\bi{2k}k}{k4^k}-\f13\sum_{k=1}^{p-1}\f{\bi{2k}k}k
\\\eq&\f23q_p(2)\pmod p
\endalign$$
since
$$\sum_{k=1}^{p-1}\f{\bi{2k}k}{k4^k}\eq 2q_p(2)\ (\mo\ p)\ \ \t{and}\ \ \sum_{k=1}^{p-1}\f{\bi{2k}k}k\eq0\ (\mo\ p^2)$$
by [ST1, (1.12) and (1.20)].
Thus, in view of the above, we get
$$\align\sum^{p-1}\Sb k=0\\k\not=(p-1)/2\endSb\f{(-3)^k}{2k+1}\eq&-2\l(\f p3\r)\f23q_p(2)+\f12\l(\f p3\r)\l(\f 83q_p(2)-q_p(3)\r)
\\=&-\l(\f p3\r)\f{q_p(3)}2\pmod p.
\endalign$$
It follows that
$$\align p\sum_{k=0}^{p-1}\f{(-3)^k}{2k+1}\eq&(-3)^{(p-1)/2}-\l(\f p3\r)\f{3^{p-1}-1}{2}
\\=&(-3)^{(p-1)/2}-\l(\f p3\r)\f{(-3)^{(p-1)/2}+(\f{-3}p)}2\l((-3)^{(p-1)/2}-\l(\f{-3}p\r)\r)
\\\eq&(-3)^{(p-1)/2}-\l((-3)^{(p-1)/2}-\l(\f{-3}p\r)\r)=\l(\f p3\r)\pmod{p^2}.
\endalign$$
We are done. \qed

\proclaim{Lemma 3.4} For any prime $p$, we have
$$k\bi{2k}k\sum_{r=0}^{p-1}\bi{-k}r\bi{-k-1}r\eq p\pmod{p^2}\quad\t{for all}\ k=1,\ldots,p-1.\tag3.5$$
\endproclaim
\Proof. Define
$$u_k=\sum_{r=0}^{p-1}\bi{-k}r\bi{-k-1}r\quad\t{for all}\ k\in\N.$$
Applying the Zeilberger algorithm via {\tt Mathematica 9}, we find the recurrence
$$\align &k(k+1)^2(2(2k+1)u_{k+1}-ku_k)
\\=&(p+k)(p+k-1)(2kp+p+3k^2+3k+1)\bi{-1-k}{p-1}\bi{-k}{p-1}
\\=&p^2\bi{p+k}p\bi{p+k-1}{p}(2kp+p+3k^2+3k+1).
\endalign$$
Thus, for each $k=1,\ldots,p-2$, we have
 $$2(2k+1)u_{k+1}\eq ku_k\pmod{p^2}$$
and hence
$$\align(k+1)\bi{2(k+1)}{k+1}u_{k+1}=&2(k+1)\bi{2k+1}{k+1}u_{k+1}
\\=&2(2k+1)\bi{2k}ku_{k+1}\eq k\bi{2k}ku_k\pmod{p^2}.
\endalign$$
So it remains to prove $\bi 21u_1\eq p\pmod{p^2}$. With the help of the Chu-Vandermonde identity, we actually have
$$\align u_1=&\sum_{r=0}^{p-1}(-1)^r\bi{-2}r=(-1)^{p-1}\sum_{r=0}^{p-1}\bi{-1}{p-1-r}\bi{-2}r
\\=&(-1)^{p-1}\bi{-3}{p-1}=\bi{p+1}{p-1}=\f{p^2+p}2.
\endalign$$
This concludes the proof. \qed

\medskip\noindent{\it Proof of Theorem 1.1}. (i) By [S12, (2.13)],
$$\f1p\sum_{k=0}^{p-1}(2k+1)A_k(x)\eq\sum_{k=0}^{p-1}\f p{2k+1}\l(1-2p^2H_k^{(2)}\r)x^k\pmod{p^4}.$$
Combining this with (3.1) we immediately get (1.4).

By [S12, (1.6)-(1.7)],
$$\f1p\sum_{k=0}^{p-1}(2k+1)A_k\eq1+\f 76p^3B_{p-3}\pmod{p^4}$$
and
$$\f1p\sum_{k=0}^{p-1}(2k+1)A_k(-1)\eq\l(\f{-1}p\r)-p^2E_{p-3}\pmod{p^3}.$$
Combining this with (3.1) we obtain (1.5) and (1.6).
In view of (1.4) and (3.4), we get (1.7).

(ii) With the help of (2.7),
$$\align\sum_{l=1}^{p-1}\f{g_l(x)}l=&\sum_{l=1}^{p-1}\f1l\sum_{k=0}^l\bi lkf_k(x)=H_{p-1}+\sum_{l=1}^{p-1}\sum_{k=1}^l\f{f_k(x)}l\bi lk
\\\eq&\sum_{k=1}^{p-1}\f{f_k(x)}k\sum_{l=k}^{p-1}\bi{l-1}{k-1}=\sum_{k=1}^{p-1}\f{f_k(x)}k\bi{p-1}k
\\\eq&\sum_{k=1}^{p-1}\f{(-1)^k}kf_k(x)(1-pH_k)\pmod{p^2}.
\endalign$$
In view of [S13b, (2.7)], this implies that
$$\sum_{k=1}^{p-1}\f{g_k(x)}k\eq p\sum_{k=(p+1)/2}^{p-1}\f{x^k}{k^2}-p\sum_{k=1}^{p-1}\f{(-1)^kH_k}kf_k(x)\pmod{p^2}.\tag3.6$$
So (1.8) follows.

By induction, for any integers $m>k\gs0$, we have
$$\sum_{n=k}^{m-1}(2n+1)\bi{n+k}{2k}=\f{m(m-k)}{k+1}\bi{m+k}{2k}.$$
This, together with (2.8) and (3.2), yields
$$\align \sum_{n=0}^{p-1}(-1)^n(2n+1)A_n=&\sum_{n=0}^{p-1}(2n+1)\sum_{k=0}^n\bi{n+k}{2k}\bi{2k}k(-1)^{k}g_k
\\=&\sum_{k=0}^{p-1}\bi{2k}k(-1)^kg_k\sum_{n=k}^{p-1}(2n+1)\bi{n+k}{2k}
\\=&\sum_{k=0}^{p-1}\bi{2k}k(-1)^kg_k\f{p(p-k)}{k+1}\bi{p+k}{2k}
\\=&g_{p-1}\bi{2p-2}{p-1}(2p-1)+p^2\sum_{k=0}^{p-2}\bi{p-1}k\bi{p+k}k(-1)^k\f{g_k}{k+1}
\\=& p\,g_{p-1}\bi{2p-1}{p-1}+p^2\sum_{k=1}^{p-1}\f{g_{k-1}}{k}
\\\eq& p\,g_{p-1}+p^2\sum_{k=1}^{p-1}\f{g_{k-1}}{k}\pmod{p^4}
\endalign$$
since $\bi{2p-1}{p-1}\eq1\ (\mo\ p^3)$ by Wolstenholme's theorem.
Combining this with (2.18) and (3.3), we obtain
$$p\l(\f p3\r)\eq p\l(\f p3\r)(1+2p\,q_p(3))+p^2\sum_{k=1}^{p-1}\f{g_{k-1}}k\pmod{p^3}$$
and hence (1.9) follows.

(1.10) follows from a combination of (1.5) and (1.11) in the case $n=p$.
If we let $u_n$ denote the left-hand side or the right-hand side of (1.11), then by applying the Zeilberger
algorithm via {\tt Mathematica 9} we get the recurrence relation
$$\align &(n+2)(n+3)^2(2n+3)u_{n+3}
\\=&(n+2)(22n^3+121n^2+211n+120)u_{n+2}
\\&-(n+1)(38n^3+171n^2+229n+102)u_{n+1}+9n^2(n+1)(2n+5)u_n
\endalign$$
for $n=1,2,3,\ldots$. Thus (1.11) can be proved by induction.

(iii) Now we show (1.12)-(1.14) provided $p>5$.

 Observe that
$$\align \sum_{l=1}^{p-1}\f{g_l(x)-1}{l^2}=&\sum_{l=1}^{p-1}\f1{l^2}\sum_{k=1}^l\bi lk^2\bi{2k}kx^k
=\sum_{k=1}^{p-1}\f{\bi{2k}k}{k^2}x^k\sum_{l=k}^{p-1}\bi{l-1}{k-1}^2
\\=&\sum_{k=1}^{p-1}\f{\bi{2k}k}{k^2}x^k\sum_{j=0}^{p-1-k}\bi{k+j-1}{j}^2=\sum_{k=1}^{p-1}\f{\bi{2k}k}{k^2}x^k\sum_{j=0}^{p-1-k}\bi{-k}j^2
\\\eq&\sum_{k=1}^{p-1}\f{\bi{2k}k}{k^2}x^k\sum_{j=0}^{p-1-k}\bi{p-k}j^2\pmod{p}.
\endalign$$
Recall that $H_{p-1}^{(2)}\eq0\ (\mo\ p)$. Also, for any $k=1,\ldots,p-1$ we have
$$\sum_{j=0}^{p-1-k}\bi{p-k}j^2=\sum_{j=0}^{p-k}\bi{p-k}j\bi{p-k}{p-k-j}-1=\bi{2(p-k)}{p-k}-1$$
by the Chu-Vandermonde identity. Thus
$$\sum_{k=1}^{p-1}\f{g_k(x)}{k^2}\eq\sum_{k=1}^{p-1}\f{\bi{2k}k}{k^2}x^k\l(\bi{2(p-k)}{p-k}-1\r)
\eq-\sum_{k=1}^{p-1}\f{\bi{2k}k}{k^2}x^k\pmod{p}$$
(Note that $\bi{2k}k\bi{2(p-k)}{p-k}\eq0\ (\mo\ p)$ for $k=1,\ldots,p-1$.) It is known that
$$\sum_{k=1}^{p-1}\f{(-1)^k}{k^2}\bi{2k}k\eq0\pmod p\tag3.7$$
(cf. Tauraso [T]) and moreover
$$\sum_{k=1}^{(p-1)/2}\f{(-1)^k}{k^2}\bi{2k}k\eq\f{56}{15}pB_{p-3}\pmod{p^2}$$
by Sun [S14]. So (1.12) is valid.
\medskip

Note that
$$\align\sum_{l=1}^{p-1}\f{g_l(x)-1}l=&\sum_{l=1}^{p-1}\f1l\sum_{k=1}^l\bi lk^2\bi{2k}kx^k
=\sum_{k=1}^{p-1}\bi{2k}kx^k\sum_{l=k}^{p-1}\f1k\bi{l-1}{k-1}\bi lk
\\=&\sum_{k=1}^{p-1}\f{\bi{2k}k}kx^k\sum_{j=0}^{p-1-k}\bi{k+j-1}j\bi{k+j}j.
\endalign$$
For $1\ls k\ls p-1$ and $p-k<j\ls p-1$, clearly
$$\bi{k+j-1}j\bi{k+j}j=\f{(k+j-1)!(k+j)!}{(k-1)!k!(j!)^2}\eq0\pmod{p^2}.$$
If $j=p-k$ with $1\ls k\ls p-1$, then
$$\align\bi{k+j-1}j\bi{k+j}j=&\bi{p-1}j\bi pj=\f pj\bi{p-1}{j-1}\bi {p-1}j
\\\eq&-\f pj\eq\f pk\pmod{p^2}.
\endalign$$
Recall that $H_{p-1}\eq0\pmod{p^2}$. So we have
$$\align\sum_{k=1}^{p-1}\f{g_k(x)}k\eq&\sum_{k=1}^{p-1}\f{\bi{2k}k}kx^k\(\sum_{j=0}^{p-1}\bi{k+j-1}j\bi{k+j}j-\f pk\)
\\=&\sum_{k=1}^{p-1}\f{\bi{2k}k}kx^k\sum_{j=0}^{p-1}\bi{-k}j\bi{-k-1}j-p\sum_{k=1}^{p-1}\f{\bi{2k}k}{k^2}x^k
\\\eq&\sum_{k=1}^{p-1}\f{x^k}{k^2}p-p\sum_{k=1}^{p-1}\f{\bi{2k}k}{k^2}x^k=p\sum_{k=1}^{p-1}\f{1-\bi{2k}k}{k^2}x^k\pmod{p^2}
\endalign$$
with the help of (3.5). Thus, in view of (3.7) we get
$$\sum_{k=1}^{p-1}\f{g_k(-1)}k\eq p\sum_{k=1}^{p-1}\f{(-1)^k}{k^2}=p\sum_{k=1}^{(p-1)/2}\l(\f{(-1)^k}{k^2}+\f{(-1)^{p-k}}{(p-k)^2}\r)\eq0\ (\mo\ p^2).$$
This proves (1.13). Combining this with (3.6) we obtain
$$\align\sum_{k=1}^{p-1}\f{(-1)^kf_k(-1)}kH_k\eq&\sum_{k=(p+1)/2}^{p-1}\f{(-1)^k}{k^2}\eq-\sum_{j=1}^{(p-1)/2}\f{(-1)^j}{j^2}
\\\eq&-2\l(\f{-1}p\r)E_{p-3}\pmod p
\endalign$$
with the help of [S11b, Lemma 2.4]. So (1.14) holds.

In view of the above, we have completed the proof of Theorem 1.1. \qed

\heading{4. Some open conjectural congruences}\endheading

In this section we pose several related conjectural congruences.

\proclaim{Conjecture 4.1} {\rm (i)} For any integer $n>1$, we have
$$\sum_{k=0}^{n-1}(9k^2+5k)(-1)^kf_k\eq0\pmod{(n-1)n^2}$$
Also, for each odd prime $p$ we have
$$\sum_{k=0}^{p-1}(9k^2+5k)(-1)^kf_k\eq3p^2(p-1)-16p^3q_p(2)\pmod{p^4}.$$

{\rm (ii)} For every $n=1,2,3,\ldots$, we have
$$\f1n\sum_{k=0}^{n-1}(4k+3)g_k(x)\in\Z[x]$$
and the number
$$\f1{n^2}\sum_{k=0}^{n-1}(8k^2+12k+5)g_k(-1)$$
is always an odd integer. Also, for any prime $p$ we have
$$\sum_{k=0}^{p-1}(8k^2+12k+5)g_k(-1)\eq3p^2\pmod{p^3}.$$
\endproclaim

For any nonzero integer $m$, the $3$-adic valuation $\nu_3(m)$ of $m$ is the largest $a\in\N$ with $3^a\mid m$. For convenience, we also set $\nu_3(0)=+\infty$.

\proclaim{Conjecture 4.2} Let $n$ be any positive integer.
Then
$$\nu_3\(\sum_{k=0}^{n-1}(2k+1)(-1)^kA_k\)=3\nu_3(n)\ls \nu_3\(\sum_{k=0}^{n-1}(2k+1)^3(-1)^kA_k\).$$
If $n$ is a positive multiple of $3$, then
$$\nu_3\(\sum_{k=0}^{n-1}(2k+1)^3(-1)^kA_k\)=3\nu_3(n)+2.$$
\endproclaim

\proclaim{Conjecture 4.3} For $n\in\N$ define
$$F_n:=\sum_{k=0}^n\bi nk^3(-8)^k\quad\t{and}\quad G_n:=\sum_{k=0}^n\bi nk^2(6k+1)C_k.$$
For any $n\in\Z^+$, the number
$$\f1n\sum_{k=0}^{n-1}(6k+5)(-1)^kF_k$$ is always an odd integer.
Also, for any prime $p>3$ we have
$$\sum_{k=0}^{p-1}(-1)^kF_k\eq\l(\f p3\r)\ (\mo\ p^2)\ \ \t{and}\ \ \sum_{k=1}^{p-1}G_k\eq-\f43p^3B_{p-3}\ (\mo\ p^4).$$
\endproclaim
\Remark\ 4.1. For any prime $p>3$, the author [S13b, S12] proved that $\sum_{k=0}^{p-1}(-1)^kf_k\eq(\f p3)\pmod{p^2}$ and
$\sum_{k=1}^{p-1}h_k\eq0\pmod{p^2}$ with $h_k=\sum_{j=0}^k\bi kj^2C_j$.

\medskip

\Ack. The author would like to thank the referee for helpful comments.

\enddocument